\documentclass[12pt]{amsart}
\usepackage{amssymb}
\usepackage{epsf,latexsym}
\newtheorem{thm}{Theorem}
\newtheorem{prop}[thm]{Proposition}
\newtheorem{cor}[thm]{Corollary}
\newtheorem{remark}[thm]{Remark}

\newenvironment{rem}{\begin{remark}\rm}{\end{remark}}

\begin{document}
\title{A special case of $\boldsymbol{sl(n)}$-fusion
coefficients}
\author{Geanina Tudose}
\address{Department of Mathematics and Statistics\\
York University\\
North York, Ont., M3J 1P3\\
CANADA}
\email{gtudose@mathstat.yorku.ca}
\begin{abstract}
We give a combinatorial description of $sl(n)$-fusion coefficients in the case where one partition has at most two
columns.
As a result we establish some properties for this case 
including solving the conjecture
that fusion coefficients are increasing with respect to the level $k$.
\end{abstract}
\maketitle

\section{Introduction}
Fusion coefficients first appeared in the literature as the structure
constants of the Verlinde (fusion) algebra associated to an affine
Kac-Moody algebra $\mathfrak{\hat{g}}$ in the Wess-Zumino-Witten model
of conformal field theory.
Since then, many equivalent interpretations have been found in 
other contexts such as quantum groups and Hecke algebras at root of unity~\cite{GdWz}, quantum
cohomology of the grassmannian~\cite{BCF}, spaces of generalized theta functions, spaces of intertwiners in vertex operator algebras~\cite{Was}, knot
invariants for 3-manifolds~\cite{T} and others.
 
 If
$\mathfrak{g}$
is a semi-simple finite
dimensional Lie algebra and $L(\lambda)$ is the integrable representation
of $\mathfrak{g}$ with
highest weight $\lambda$, the tensor product coefficients $N_{\lambda \mu}^{\nu}$ are
defined by the relation $L(\lambda )\otimes L(\mu )=\oplus N_{\lambda \mu}^{\nu}L(\nu)$.
For $\mathfrak{g}=sl(n)$ the integrable representations are indexed by
partitions and the tensor product coefficients are the well-known
Littlewood-Richardson coefficients
$c_{\lambda \mu}^\nu$. Given a positive level $k$, the fusion
coefficients
$N_{\lambda
\mu}^{(k)
\nu}$ are defined by 
$$
L(\lambda ) \otimes_k L(\mu )=\oplus N_{\lambda \mu}^{(k) \nu}L(\nu)
$$
where the fusion product $\otimes_k$ is the reduction of the tensor
product via the representation at level $k$ of the algebra $\mathfrak{\hat{g}}$. A more
detailed approach to fusion coefficients arising in conformal field
theory is given in~\cite{W1,W2}. For our purposes we will give in
Section 3 an equivalent definition for the case $\mathfrak{g}=sl(n)$.

By some representation theoretic arguments it is known that these
coefficients are non-negative but a general combinatorial description is
still lacking even for type $A$. Only some particular cases are known: the cases $n=2$ 
and $n=3$~\cite{BMW,BKMW,KMSW} where the combinatorial objects used
are the Berenstein-Zelevinski triangles, and more recently the case
where all partitions in the product are rectangles~\cite{SS1, SS2} in
which affine crystal theory for perfect crystals was used. In addition,
a $q$-analogue of fusion coefficients has also been introduced~\cite{FLOT}. 

To date, the most effective algorithm for computing fusion coefficients
for any type is the Kac-Walton algorithm ~\cite{Kac,W1}. In this
algorithm, the fusion coefficients are expressed
 in terms of the
tensor
product coefficients
$$
N_{\lambda \mu}^{(k) \nu}=\ \sum_{\stackrel{\mbox{\scriptsize $w \in \hat{W}_k$}}%
{w. \nu \in P^{+}}}det(w) N_{\lambda\mu}^{w. \nu}
$$
where $\hat{W_k}$ is the affine Weyl group, $w. \nu=w(\nu+\rho)-\rho$, $P^+$ is the set of dominant weights
and $\rho$ is the sum of fundamental weights. In this notation the affine Weyl
groups
are isomorphic, and only the action of the reflection $s_0$ on
the weight lattice of the algebra $\mathfrak{\hat{g}}$ is different with respect
to the level
$k$ i.e.
$$
 s_0(\hat{\lambda})=\hat{\lambda}-(k-(\lambda,\theta))\hat{\alpha}_0  
$$
where $\theta$ is the highest root of $\mathfrak{g}$, $\{ \hat{\alpha}_i,\
i=0,1,\ldots n-1\}$ are the simple roots, and $(\cdot|\cdot)$ is the
symmetric billinear form on the Cartan subalgebra of $\mathfrak{\hat {g}}$.  

In this paper we use the
interpretation given by Goodman and Wenzl~\cite{GdWz} to give a combinatorial description for $sl(n)$ fusion
coefficients where the partition $\mu$ has two columns.

  Our main result is Theorem 12 where we show that the
  fusion coefficients count paths in the Young's
  lattice with some extra conditions. An equivalent interpretation in
  terms of Littlewood-Richardson tableaux is given in Remark 13. The tool for finding this description is the pairing technique for proving the classical Littlewood-Richardson rule by means of a sign-reversing involution.
Therefore we include in Section 2 a proof of the classical rule
so that in Section 3 we can construct the involution for fusion coefficients.
In Section 4 we establish some interesting properties of these
coefficients in our specific case. Some of these confirm known
  properties such as positivity and the inequality $N_{\lambda \mu}^{(k) \nu}\leq
N_{\lambda\mu}^{\nu}$, but most importantly we confirm the increasing property as function of $k$
conjectured in~\cite{W2} i.e. $N_{\lambda \mu}^{(k) \nu} \leq N_{\lambda
\mu}^{(k+1) \nu}$. We conclude our paper with Section 5 where we
  propose another avenue for approaching the problem.

\section { Proof of the Littlewood-Richardson rule}

 The proof of the LR-rule is based on the Jacobi-Trudi determinantal
identities and uses a sign reversing involution which yields a
combinatorial characterization of the LR-coefficients in terms of paths in
the Young's lattice. The involution is an adaptation of the involution
constructed by Remmel and Shimozono~\cite{RS}.

  The LR-coefficients are the structure constants $ c_{\lambda \mu}^{\nu} 
$
for the ring of symmetric polynomials with respect to the basis of Schur
functions:
\begin{equation}\label{eq:one}
         s_{\lambda} s_\mu =\ \sum_\nu c_{\lambda\mu}^\nu s_\nu.
\end{equation}

We intend to give a characterization of these coefficients of the form
 $$
       c_{\lambda\mu}^\nu =\ \sharp\left\{ \begin{array}{ll}
                               \hbox {paths in the Young's lattice from }
\lambda\hbox { to } \nu\\
                         \hbox { satisfying conditions imposed by } \mu
 \end{array}
 \right\}.
 $$

 There are many ways of getting to this result depending on which
determinantal formula we use. We shall choose the one expressing  the Schur
functions in terms of the elementary symmetric polynomials $ e_k $.
The reason for this choice is accounted for in the proof of the rule
for fusion coefficients.

 In order to prove the LR-rule we first need some definitions. Most of
 those not given here and results concerning symmetric functions that
 we use can be found in~\cite{M}. For a partition $\lambda$ we consider
its diagram to be the set of
points $(i,j) \in {\mathbb{Z}}^2$ such that  $1\leq i\leq \lambda_j$, 
where $1 \leq j \leq length(\lambda)$.

 We denote {\sl a path} $P$ in the Young's lattice from $ \lambda $ to
$\nu$  by
{\sl a chain of partitions}
$$
    P\quad :\quad \lambda^{(0)}=\lambda \subseteq \lambda^{(1)} \subseteq
\cdots
\lambda^{(n)}=\nu
$$
 where each partition $ \lambda^{(k)} $ differs from the previous one
$ \lambda^{(k-1)} $ by exactly one box. We also denote by $|P|=n$ the
length of the path $P$.

  Sometimes we need paths from $\lambda$ to $\nu$ made from {\sl successive
paths} i.e.\ $P=\ P_1*P_2*\cdots *P_m$, where each $P_i$ is a path from 
$\lambda^{(i)}$ to $\lambda^{(j)}$ with $i\leq j$, and $*$ denotes the
concatenation of the paths.

$$
\lower0.2in\hbox{$\epsfxsize=0.45in\epsfbox{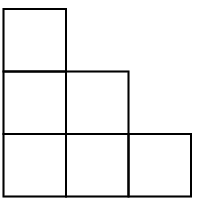}$}\
\ \hookrightarrow \ \ 
\lower0.2in\hbox{$\epsfxsize=0.6in\epsfbox{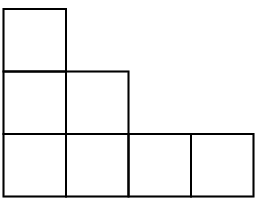}$}\ \ 
\hookrightarrow\ \ 
\lower0.2in\hbox{$\epsfxsize=0.6in\epsfbox{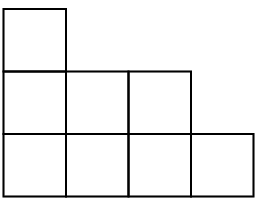}$}\
\ \hookrightarrow \ \ 
\lower0.2in \hbox{$\epsfxsize=0.6in\epsfbox{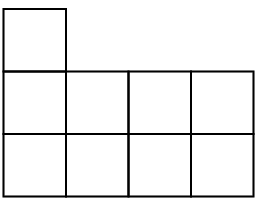}$}$$
$$
\lambda^{(0)}\ \ \ \ \ \hookrightarrow \ \ \ \ \ 
\lambda^{(1)}\ \ \ \ \
\hookrightarrow \ \ \ \ \  \lambda^{(2)}\ \ \ \ \
\hookrightarrow\ \ \ \ \ \ \ \lambda^{(3)}
$$

\centerline{{\bf Figure 1.  } $P:\ \lambda^{(0)}\subseteq
\lambda^{(1)}\subseteq\lambda^{(2)}\subseteq\lambda^{(3)}$
{\sl and  we can also write, say }}
$$P=P_1*P_2,\ \hbox{\sl where }
P_1=\lambda^{(0)}\subseteq\lambda^{(1)}\subseteq\lambda^{(2)}\ \hbox{\sl
and } P_2=\lambda^{(2)}\subseteq \lambda^{(3)}. $$
\vskip 12pt
  Next we introduce a {\sl labeling} of each box in a partition in order
to define a 1-1 correspondence between the paths and the sequence of
labels such that the boxes on the diagonals $x-y=i$ are indexed by $i$.

$$
\epsfxsize=1in\epsfbox{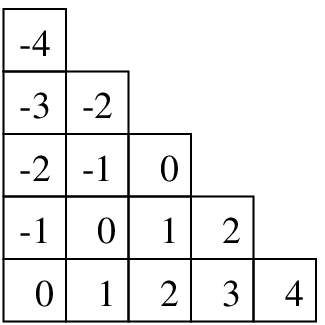}
$$
\centerline{{\bf Figure 2.} {\sl Labeling of $\lambda=(5,4,3,2,1)$.}} 
\vskip 12pt
\noindent
Using this labeling we identify the path
$P=\lambda\subseteq\lambda^{(1)}\subseteq\cdots \lambda^{(n)}$ with the
sequence
of boxes added in each step. From here  we
shall write the labels of these boxes as
$$
l(P)=\ (l_1,\,l_2,\ldots l_n),
$$
where $l_i$ is the label of $\lambda^{(i)}/\lambda^{(i-1)}$, for
$i=1,\ldots,n$.
We say that $P$ is {\sl a decreasing path } if $l(P)$ is decreasing. If
$\alpha$ is a
sequence of integers $(\alpha_1,\alpha_2,\ldots \alpha_k)$ with
$\alpha_1+\alpha_2+\cdots \alpha_k=|P|$ 
we say
that {\sl $P$ has ascents included in positions $\alpha$} if
$$(l_{\alpha_0+\alpha_1+\cdots
+\alpha_i+1},\, l_{\alpha_0+\alpha_1+\cdots+\alpha_i+2}, \ldots,
l_{\alpha_0+\alpha_1+\cdots +\alpha_{i+1}})$$
 is a
decreasing
sequence for every $i \in
\{0,1,2,\ldots ,k-1\}$, where $l_j$ is the $j^{th}$ component of $l(P)$ and $\alpha_0=0$.

\noindent We make the
convention that if $\alpha$ contains negative integers the set of paths
with ascents included in positions $\alpha$ is the empty set. We also note
that a path $P$ can have ascents included in different $\alpha$'s.
 
{\sl Example 1: $P$ in {\bf Figure 1} has  $l(P)=\ (3,1, 2)$ and $P$ has
ascents
in
$\alpha=(2,1)$ and also \indent ascents included in positions $(1,1,1)$}.

\noindent
 For each general path $P$ and a sequence $\alpha=(\alpha_1,\ldots
\alpha_k)$ such that $P$ has ascents included in positions $\alpha$, we
cut
the
path $P$ into $k$ consecutive paths each of length $\alpha_i$, $i=1
\ldots k$; then we associate a tableau $T_P$ whose  columns $i$ are made
from the sequence
of labels of $P_i$ written top-to-bottom. Sometimes when $\alpha$ is
understood we will make no distinction between $P$ and $T_P$.
 We say that {\sl a path $P$ fits a partition $\mu$}, if
$T_P\in CS(\mu)$,
where $CS(\mu)$ represents the Young tableaux of shape $\mu$, strictly
increasing in columns and weakly increasing in rows.

\begin{thm}
(Littlewood-Richardson rule)
\begin{equation}\label{eq:two}
        c_{\lambda\mu}^\nu=\ \sharp\lbrace\hbox{paths $P$ from $\lambda$
to $\nu$ that fit $\mu$}\rbrace.
\end{equation}
\end{thm}
\noindent{\bf Proof. }

Let $\mu'$ denote the conjugate partition of $\mu$.

Using the Jacobi-Trudi identity to express $s_\mu$ in terms of the elementary symmetric
 functions given in Equation~(\ref{eq:one}) we
get 
$$
   s_{\lambda}s_{\mu}=\  s_{\lambda} det(e_{\mu'_i-i+j})_{1 \leq i,j \leq n}=\ \sum_\nu
c_{\lambda
\mu}^\nu s_\nu \qquad\hbox{where $n\geq l(\mu')$}
$$
and when we expand the determinant we have
\begin{equation}\label{eq:three}
s_\lambda \sum_{\sigma \in {\mathcal {S}}_n}(-1)^\sigma e_{\sigma. \mu'}=\ \sum_\nu c_{\lambda
\mu}^{\nu} s_\nu
\end{equation}
where $\sigma. \mu'=\ \sigma (\rho+\mu')-\rho$ and $\rho=\ (n-1,n-2,\ldots,1,0).$
 On the other hand multiplying a Schur function with an elementary symmetric
function we get
$$
    s_{\lambda} e_k=\ \sum_{\nu /\lambda=k-column\, strip} s_\nu.
$$
We can also view this equality in terms of paths in the Young's lattice
as
\begin{equation}\label{eq:four}
  s_{\lambda} e_k =\ \sum_\nu a_{\lambda(k)}^\nu s_\nu
\end{equation}
 where $ a_{\lambda(k)}^\nu=\ \sharp
\lbrace \hbox{decreasing paths from $\lambda$ to $\nu$ of length $k$} \rbrace. $
It is not difficult to see that, indeed
\[ a_{\lambda(k)}^\nu= \left\{ \begin{array}{ll}
                             1 & \hbox { if
                              $\nu$/$\lambda$ is a $k$ column strip }\\
                             0 & \hbox{ otherwise}.
                           \end{array}
                           \right. \]

Using rule~(\ref{eq:four}) repeatedly, the left-hand side
of Equation~(\ref{eq:three}) becomes
$$
 \sum_{\sigma\in {\mathcal{S}}_n}(-1)^\sigma \sum_\nu s_\nu e_{\sigma. \mu '}=
\ \sum_{\sigma \in{\mathcal{S}}_n}\sum_\nu (-1)^\sigma a_{\lambda
(\sigma. \mu')}^\nu s_\nu
$$
where
$$
 a_{\lambda (\sigma. \mu ')}^\nu =\ \sharp\{\hbox{paths from
$\lambda$ to $\nu$ with ascents in positions $\sigma. \mu'\}$.}
$$

 Hence 
\begin{equation}\label{eq:five}
  c_{\lambda\mu}^\nu = \ \sum_{\sigma \in {\mathcal{S}}_n}(-1)^\sigma
a_{\lambda(\sigma. \mu')}^\nu =\ \sum_{(\sigma,P)\in \Omega}(-1)^\sigma
\end{equation}
where $\Omega$ is the set of pairs $(\sigma,P)$, $ \sigma\in {\mathcal{S}}_n$ and  
$P$ is a path from $\lambda$ to $\nu$ with ascents in
positions $\sigma. \mu'$. Since $P$ has ascents included in
positions
$\sigma. \mu '$, we can write $P=\ P_1 *\cdots * P_n$  in which $P_i$
is a decreasing path of length $|P_i|=\ (\sigma. \mu ')_i$.

The next step is to construct a sign reversing involution on the set $\Omega$. The involution uses the
crystal operators defined by a {\sl pairing}
and was constructed in~\cite{RS} which also contains further details.

Suppose that a path $P$ is made from two successive paths 
$P=P_1* P_2$ of lengths $p$ and $q$.
The {\sl word} of $P$ denoted by $w$ is the sequence of
all labels in $P$ sorted in increasing order. A label can appear at most
twice, i.e.\ once in each column and if this
happens we consider the first occurrence corresponding to the first
column in $T_P$ and the second occurrence corresponding to
the second column. 

We construct $\hat{w}$ in the following way

--- Replace every letter in $w$ which is a label in $l(P_1)$ by a  left parenthesis

--- Replace every letter in $w$ which is a label in $l(P_2)$ by a right parenthesis.

{\sl Example 2: The path $P$ with $T_P$ =  \lower0.25in\hbox{$\epsfxsize=0.4in
\epsfbox{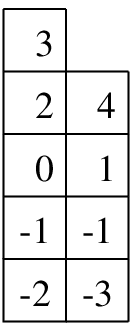}$} has the word
$w=\bar{3}\,\bar{2}\,\bar{1}\,\bar{1}\,0\,1\,2\,3\,4$, where $\bar{n}=-n$.
\indent The parentheses structure is
 $$\hat{w}=)(()()(().$$}

\noindent
We say that a letter is {\sl paired } if it corresponds to a parenthesis
that
is matched under the
usual rule of parenthesization. Otherwise we call it {\sl unpaired}. We
say that a word $w$ has $
type\,(l,r)$  if there are
$l$ unpaired left parenthesis and $r$ unpaired right parenthesis.

Next we define two operators on words which will be partial functions, 
the raising
operator $e$ and the lowering operator $f$ where 

--$e$ changes the rightmost unpaired right parenthesis into a left one.

--$f$ changes the leftmost unpaired left parenthesis into a right one.

\noindent It is clear that for $e$ or $f$ to be applied we need $r>0$
(resp. $l>0$).
We shall also write $e(P)$ or $f(P)$ and understand that $e$ or $f$ is applied to the word of $P$ with
$e(P)=P_1'*P_2'$, where $|P_1'|=|P_1|+1,\ |P_2'|=|P_2|-1$ and also
$f(P)=P_1''*P_2''$ with $|P_1''|=|P_1|-1,\ |P_2''|=|P_2|+1$.

{\sl Example 3: For $P$ in Example 2, $e(w),\ f(w)$ have the parentheses
structure
$$
\widehat{e(w)}=((()()(()\ \  and\ \  \widehat{f(w)}=))()()(()
$$
\indent and the results of these operators on the tableau of $P$ are 
$$T_{e(P)}=\ \ \lower0.25in\hbox{$\epsfxsize=0.4in
\epsfbox{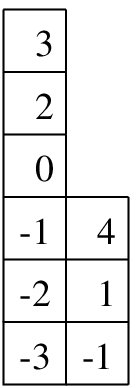}$}\ \ \hbox{and}\ \
T_{f(P)}=\lower0.25in\hbox{$\epsfxsize=0.4in
\epsfbox{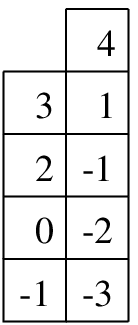}$}\ .
$$
}

The next result helps us to establish that $e$ and $f$ define an
involution.
\begin{prop}[Proposition~3 of~\cite{RS}]
 If $\eta$ is any of the operators $e,\, f$, then the unpaired subwords of $\eta (w)$ and of $w$ occupy
the same positions (assuming $\eta$ is defined) and if $w$ has at least $m$ unpaired left
parentheses then $f^me^m(w)=w$. A similar property holds for the unpaired right parentheses.
\end{prop}
Therefore we can consider that $e^{-1}=f$ and $f^{-1}=e$ where they
are defined. The following useful result is a reformulation of
Proposition~5 of~\cite{RS}.
\begin{prop}
A path $P$ fits $\mu$ ( i.e.\ $T_P\in CS(\mu)$ ) if and only if there
are no unpaired right parentheses for every two columns $(P_i,
P_{i+1})$ in $\mu$, where $i=1,\ldots \mu_1-1$.
\end{prop}
\begin{rem}
If $P$ does not fit a partition then there exists two consecutive columns
$P_i$ and $P_{i+1}$ for which we have

--at least $|P_{i+1}|-|P_i|-1$ unpaired right parentheses, if
$|P_{i+1}|-|P_i|-1>0$

--at least $-(|P_{i+1}|-|P_i|-1)$ unpaired left parentheses, if
$|P_{i+1}|-|P_i|-1<0$.

\noindent
It is an easy consequence of the expansion of the determinant
in~(\ref{eq:three}) that $|P_{i+1}|-|P_i|-1 \neq 0$.
\end{rem}
\noindent

We use Proposition~3 to construct the involution $\Psi$ on the right-hand
side
of~(\ref{eq:five}) as follows.
\begin{enumerate}
\item If  $T_P\in CS(\mu)$  then $\sigma=id$ and define 
$\Psi(id,P)=\ (id,P)$.
\item If  $T_P\notin CS(\mu )$  then let $(r,r+1)$ be the pair of
consecutive columns where a
violation of the column-strict tableau property
occurs while reading $T_P$ from right to left, bottom to top,
row-wise. We call this position canonical.
Define 
$$
\Psi(\sigma,P) :=
((r,r+1). \sigma,P_1*\cdots*P_{r-1}*e^{|P_{r+1}|-|P_r|-1}(P_r*P_{r+1})*\cdots*P_m).
$$

\end{enumerate}

\noindent
 We must show that $\Psi$ is a well-defined involution. We first check
 that $\Psi(\sigma,\,P)\in \Omega$.
This is trivial when  $T_P\in CS(\mu)$, so we shall assume that $T_P\notin CS(\mu)$. From Proposition~3
and the Remark 4, it is clear that we can define the operator
$e^{|P_{r+1}|-|P_r|-1}(P_r*P_{r+1})$.
If $\Psi(\sigma,P)=\ (\sigma ',P')$, then $P'=\ P_1*\cdots *P'_r * P'_{r+1}*\cdots *P_m$ is 
indeed a path with ascents included in positions $\sigma
'. \mu '$, since both $P'_r$ and
$P'_{r+1}$ are decreasing paths and $|P'_i|=\ (\sigma '. \mu ')_i$, for
any $i$. Thus $\Psi$ is well-defined.

 Next we shall show that $\Psi$ is an
involution. This is again obvious for $P$ a partition that fits $\mu$. Let
$P$
be a path such that $T_P\notin CS(\mu)$ and let
$(\sigma',P')=\Psi(\sigma,P)$. To see that $\Psi(\sigma',P')=(\sigma,P)$ it is necessary to show that the violation of the
column-strict tableau property occurs in the same place for both $T_P$ and
$T_{P'}$.
This violation can be either a non-increasing pair on a row $k$ and columns
$r$ and $r+1$, or the associated tableau is not a shape.

 If $T_P$ satisfies first situation, then
all the $i^{th}$ columns, $i\ne r,\,r+1$ in $T_{\Psi(P)}$ remain
unchanged. For the columns $r$ and
$r+1$
$$\epsfxsize=0.9in \epsfysize=0.95in \epsfbox{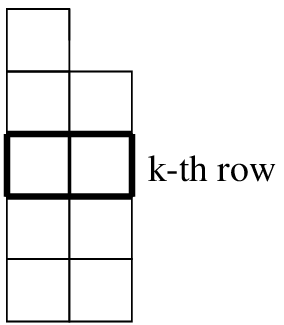}$$
\noindent
everything under the $k^{th}$ row is also unchanged since all these labels are paired
parentheses in the word of $P_r*P_{r+1}$. For this row
the only change that can occur is of the type
$$\epsfxsize=0.4in \epsfysize=0.35in \epsfbox{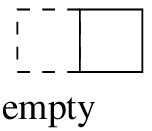}$$
\noindent
which is a violation of the shape property. Thus the canonical
position for $T_{\Psi(P)}$ is the same as for $T_P$.

If $T_P$ satisfies second situation, i.e. a violation of the shape property
occurs on the $k^{th}$ row and the columns $r$ and $r+1$, we have the reverse
of the above situation and as a result the same canonical position for
both the tableau and the image. 

Therefore $\Psi$ is a well-defined involution, is sign-reversing and
by definition, its only fixed points are $(id,P)$, where
$P\in
CS(\mu)$. This proves the characterization given in Theorem~1.
\begin{rem}
  The characterization of LR-coefficients given here is equivalent to
  the characterization where $c_{\lambda\mu}^\nu$ counts the number of
row-strict tableaux of shape $\nu/\lambda$, content $\mu'$ whose word is
lattice (read column-wise). We say that a word is lattice if every
initial subword has
(the number of $i$'s) $\geq$ (the number of $(i+1)$'s), for every $i$. 
To see the equivalence we note the 1-1 correspondence obtained by
  labeling all boxes from $P_i$ with $i$, for any $i$.
\end{rem}

  \section {The LR-rule for fusion coefficients where one partition has at
most two columns}

 The fusion coefficients we consider are the structure constants for
the fusion
algebra of WZW conformal field theories associated to
 $\widehat{sl(n)}$ at level $k$.
 This algebra $\mathcal{F}^{(n,k)}$ is isomorphic
to the algebra of symmetric polynomials $\mathbb{Q}(x_1,\ldots x_n)^{\mathcal{S}_n}/\mathcal{I}^{(n,k)}$
where
$\mathcal{I}^{(n,k)}$ is the
ideal of $\mathbb{Q}(x_1,\ldots x_n)^{\mathcal{S}_n}$ generated by the Schur functions $s_\lambda$
for which $\lambda_1-\lambda_n=k+1$, and $s_{(1^n)}-1$. The interpretation of the fusion
algebra and many results that we will use here rely on the paper of Goodman and
Wenzl~\cite{GdWz}.

Before we proceed we require some more notation and definitions most
of which can also be found in ~\cite{GdWz} or~\cite{FLOT}. In fact the
interpretation of the fusion algebra we use is taken from~\cite{GdWz} as
are many results which we will manipulate. 

We say that a partition $\lambda$ is {\sl $(n,k)$-restricted}, if  
$l(\lambda)\leq n$ and $0<\lambda_1-\lambda_n\leq k$. We denote the set of $(n,k)$-restricted
partitions by
$\Pi^{(n,k)}$. If $\lambda$ is such that $l(\lambda)\leq n$ and
$ \lambda_1-\lambda_n =\ k+1$ we call it {\sl a border diagram} and if
$\lambda$ is such that $l(\lambda)\leq n$ and $\lambda_1-\lambda_n =\ k$ 
we call it {\sl an edge diagram}.

\noindent
We say that a row-strict tableau $T$ is
{\sl $(n,k)$-restricted} if the shape of $T$ is a $(n,k)$-restricted
partition and the row-strict property is
preserved when we align the $n^{th}$ row and the first row on the
right of $k$ boxes. We denote by
$RS\Pi^{(n,k)}(\lambda,\ \mu)$ the
set of row-strict 
$(n,k)$-restricted tableaux of shape $\lambda$ and content
$\mu$. Similarly, we define the  column-strict $(n,k)$-restricted
tableaux and denote their set by $CS\Pi^{(n,k)}(\lambda,\ \mu)$.

{\sl Example: A row-strict $(4,4)$-restricted tableau}
$$
\epsfbox{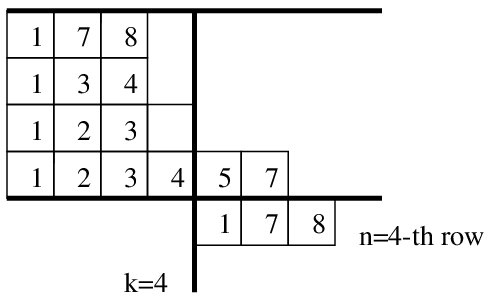}\ .
$$
The fusion algebra $\mathcal{F}^{(n,k)}$ has a linear basis indexed by the
set $\bar{\Pi}^{(n,k)}=\{\lambda,\ l(\lambda)\leq n-1,\ \lambda_1\leq
k\}$.
We can define the quotient map in the following way
$$
\Pi^{(n,k)}\ \ \longrightarrow \ \ \bar{\Pi}^{(n,k)}
$$
$$
\ \ \ \qquad\hbox{$\lambda$}\ \rightarrow\ \
\bar{\lambda}=(\lambda_1-\lambda_n,\cdots,\lambda_{n-1}-\lambda_n).
$$
\noindent
The product of two Schur functions indexed by $\bar{\Pi}^{(n,k)}$ can be recovered from the
product of Schur functions indexed by $\Pi^{(n,k)}$
. Therefore we can instead work 
with the basis $\{s_\lambda\}_{\lambda \in \Pi^{(n,k)}}$.

The structure constants of the fusion algebra are defined by
\begin{equation}\label{eq:first_def}
    s_\lambda s_\mu =\ \sum_\nu N_{\lambda\mu}^{(k)\nu}s_\nu
\qquad\hbox{where $s_\lambda,\,s_\mu,\,s_\nu\,\in \Pi^{(n,k)}$}.
\end{equation} 
 
By their equivalent interpretation to the Hecke algebras at root of
unity~\cite{GdWz}
it is known that these coefficients are nonnegative. Here we are able to
give a
combinatorial characterization for them in the case $\mu_1\leq 2$ and
in addition prove some properties one of which was conjectured
in~\cite{W2}. Using the notations from Lie algebras this means that the
weight $\mu$ has the form $\mu= \Lambda_i + \Lambda_j$, where $1\leq
i,j \leq n-1$ and $\Lambda_i$ are the fundamental weights of
$sl(n)$. 
In order to proceed we need the following result from~\cite{GdWz}.

\begin{prop}[Corollary~3.3 of~\cite{GdWz}]
If $\mu\in \Pi^{(n,k)}$, then $ s_\mu =\ det(e_{\mu_i'-i+j})_{1\leq
i,j\leq
m}$ where $m\geq l(\mu ')$ and  $e_r =\ 0$ for $r>n$ or $r<0$.
\end{prop}

 Multiplying a Schur function by an elementary symmetric function within
the fusion algebra (Proposition~2.6 of~\cite{GdWz}) we get  
\begin{equation}\label{eq:second_def}
   s_\lambda e_r =\ 
\sum_{\stackrel{\mbox{\scriptsize $\nu/\lambda =r-column\,strip$}}%
{\nu\in \Pi^{(n,k)}}} 
s_\nu.
\end{equation}

\noindent
If in Equation~(\ref{eq:first_def}) we have $\mu_1=1$ and hence
$s_\mu=e_r$, then the above expression gives the fusion coefficients
to be
$$
N_{\lambda \mu}^{(k)\nu}= \left\{ \begin{array}{ll}
                         1 & \hbox{ if $\nu$/$\lambda$ is a $r$-column
strip and $\nu \in \Pi^{(n,k)}$}\\
                         0 & \hbox{ otherwise }.
                         \end{array}
                         \right. 
$$
\noindent
For $\mu_1>1$, using  Proposition 6 on the left-hand side of Equation~(\ref{eq:first_def}) we
obtain
\begin{equation}\label{new1}
s_\lambda s_\mu =\ s_\lambda det(e_{\mu_i'-i+j})_{1\leq i,j\leq m}.
\end{equation}
By expanding the determinant and using Equation~(\ref{eq:second_def}) we get 
\begin{equation}\label{new2}
 \sum_{\sigma\in S_m}(-1)^\sigma \sum_{\nu\in\Pi^{(n,k)}}s_\nu
e_{\sigma. \mu'} =\ \sum_{\sigma\in S_m}\sum_{\nu\in
\Pi^{(n,k)}}(-1)^\sigma a_{\lambda_(\sigma. \mu')}^{(k)\nu} s_\nu,
\end{equation}
 where $a_{\lambda (\sigma. \mu')}^{(k)\nu} =\ \sharp\{$ paths
in the Young's lattice from
$\lambda$ to $\nu$ with ascents included in positions $(\sigma. \mu')$ 
and for which the partitions corresponding to these positions are
$(n,k)$-restricted$\}$.

\noindent
We denote by $\mathcal{P}^{(n,k)}_{(\sigma. \mu')}$ the set of
all such paths.

\noindent
If we equate the coefficient
of $s_\nu$ in both Equations~(\ref{eq:first_def}) and~(\ref{new2}) we get 
\begin{equation}\label{inv}
    N_{\lambda \mu}^{\nu (k)} =\ \sum_{(\sigma,P)\in \Omega_k}(-1)^\sigma
\end{equation}
 where  $\Omega_k$ is the set of pairs $(\sigma,P)$ , $\sigma\in \mathcal{S}_m$
and $P\in \mathcal{P}^{(n,k)}_{(\sigma. \mu')}$.

Our aim is to construct an involution $\Phi$ on the set $\Omega_k$
that cancels the negative terms on the right-hand side of Equation~(\ref{inv})
and that  will yield a
combinatorial description for the coefficients $N_{\lambda \mu}^{\nu (k)}$.

In this paper we consider $\mu_1=m=2$.
\begin{rem}
We exclude here the case $l(\mu)=\mu_1'=n$. In this case $s_\mu=e_ne_{\mu_2'}$,
so $N_{\lambda \mu}^{(k)\nu}=a_{\lambda
(\sigma. \mu')}^{(k)\nu}=card\mathcal\,{P}^{(n,k)}_{(n,\mu_2')}$.
\end{rem}
We may assume in what follows that $l(\mu)<n$.

\noindent
Let $\lambda$ and $\nu$ be two $(n,k)$-restricted partitions and $P$ a
decreasing path from $\lambda$ to $\nu$ with labels $l(P) =
(l_1,\ldots l_t)$ , where $l_1>l_2\ldots >l_t$, so that $\nu/\lambda$ is a
column strip.
We say that $P$ has {\sl $\bot$-label} if $P$ has a label corresponding
to the first row of the diagram $\lambda$ and $P$ has {\sl a $\top$-label}
if
there is a label corresponding to the $n^{th}$ row of the diagram $\lambda$, where $n$ is given by
the definition of $(n,k)$-restricted partition. We will denote
these labels simply by $\bot$ and $\top$.

{\sl Example 1: Suppose $\nu/\lambda$ =\ \lower0.4in\hbox{ $\epsfxsize=0.8in
\epsfbox{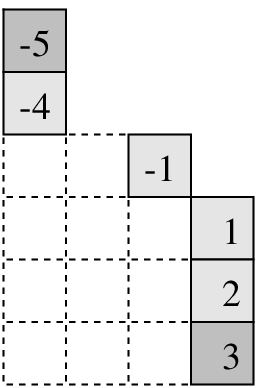}$}, where
$n=6$ and $l(P)=(3,2,1,\bar{1},\bar{4},\bar{5})$, so $3$ \indent represents the
$\bot$-label and $\bar{5}$ represents the $\top$-label}.
\vskip 12pt
\noindent
We shall write these labels in the tableau of $P$ as
\begin{equation}
\epsfxsize=0.2in \epsfbox{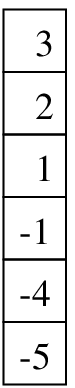}\rlap{\hskip
  3pt\lower -0.95in \hbox{$\leftarrow \bot$}}\hskip
  3pt\lower-0.07in\hbox{$ \leftarrow \top$}\ . 
\end{equation}

Since $\bot$ is the largest label and $\top$ is the
smallest, in figures where we do not specify the filling, we omit the symbols and
we just use grey boxes to indicate their presence. The four possible situations are
\begin{equation}\label{cases}
\epsfysize=0.9in \epsfbox{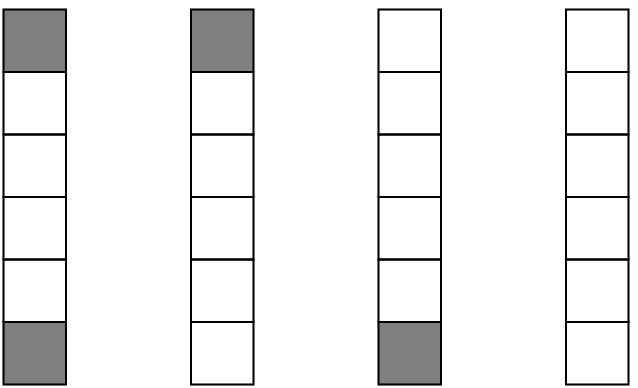}\ .
\end{equation}
The involution that we construct primarily uses the crystal
operators introduced in the previous section, however when this is not possible we define
a new modified operator.

\noindent
Let $P$ be a path $P=\ P_1*P_2 $ from $\lambda$ to $\nu$ with the
intermediate diagram $\eta$ such that all $\lambda,\ \eta,\ \nu$ are
$(n,k)$-restricted partitions, that is
\begin{equation}\label{eta}
\lambda\ \buildrel P_1 \over \longrightarrow\ \eta \
\buildrel P_2 \over \longrightarrow\
\nu.
\end{equation}

\noindent
Consider the following sets:

-- $\mathcal A = \{\hbox{paths $P=P_1*P_2 \in \mathcal{P}_{((21). \mu')}^{(n,k)}$, such that $|P_1|<|P_2|$} \}$ 

-- $\mathcal{B}= \{\hbox{paths $P=P_1*P_2 \in \mathcal{P}_{(\mu')}^{(n,k)}$, such that $|P_1|\geq
   |P_2|$} \}$. 

\noindent
It is not difficult to see that the set
$\Omega_k$ of Equation~(\ref{inv}) is in fact 
$$
\Omega_k=\{( (21),\,P),\ P \in \mathcal A \} \cup \{ ( id,\, P),\ P
\in \mathcal B\}
$$
\noindent 
and by an abuse of notation we will write $\Omega_k= \mathcal A \cup
\mathcal B$. The involution $\Phi$ that we will construct will have
the property that $\Phi(\mathcal A) \subseteq \mathcal B$ and
$\Phi(\mathcal B) \subseteq \mathcal A$. We will start be defining $\Phi$ on the set $\mathcal A$. 

\noindent
If $P \in \mathcal A$ is a path as in~(\ref{eta}) the image
$\Phi(P)=P_1'*P_2'$ will be
 $$
\lambda\ \buildrel P_1' \over \longrightarrow\ \eta' \
\buildrel P_2' \over \longrightarrow\ \nu.
$$
\noindent
Since we want $\Phi(P) \in \mathcal B$ we must ensure that $\eta'$ is
a restricted partition. We denote by $\Psi$ the involution for the
classical LR-rule constructed previously. We consider the following two cases. 

{\bf Case 1. } Suppose
$\nu$ is not an edge diagram i.e.\ $\nu_1-\nu_n<k$.

In this case let $\Phi(P)=\ \Psi(P)$. To show that $\Phi$ is well-defined recall that the rightmost
$|P_2|-|P_1|-1$ unpaired right parentheses from the word of $P$ must change into left parentheses
and hence this number of labels from the column $P_2$ move into the first column. We note that
if the largest label of $P_2$ corresponds to an unpaired parenthesis,
then this is the first to
move. We must therefore check that if this unpaired label is $\bot$ we
still obtain a partition $\eta'\in \Pi^{(n,k)}$. When the
$\Psi$-operator is
applied to $P$, the image of the intermediate partition denoted by
$\eta'$ satisfies
$\eta_1'-\eta'_n=\eta_1-\eta_n+\{-1,\,0,\,1\}$. Now since
$\nu_1-\nu_n<k$ we only
need to see what happens when $\nu_1-\nu_n=k-1$ and
$\eta_1-\eta_n=k$. In other words we have $\top \in P_2$, and $\bot
\notin P_2$ i.e.
$$ T_{P_2}=\ \ \lower 0.4in\hbox{$\epsfysize=0.9in
  \epsfbox{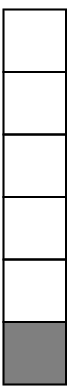}$}\ .
$$
In this case the length of the first row of the diagram $\eta'$ will be equal to the length of the
first row of the diagram $\eta$, so $\eta'\in \Pi^{(n,k)}$, too. Regardless of the presence of
$\bot$ or $\top$ in $P_1$, the image
$$ T_{\Phi(P)}=\ \ \lower 0.4in\hbox{$\epsfysize=0.75in 
\epsfbox{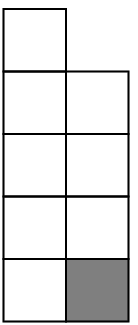}$}$$
\noindent
is not a column strict tableau because it belongs to the image of the
operator $\Psi$.

\noindent
{\bf Case 2. } The partition $\nu$ is an edge diagram i.e.\ $\nu_1-\nu_n=k$.

From~(\ref{cases}) it follows that there are $16$ cases to be studied depending on whether $\bot$ or $\top$ appears in
$P_1$ or $P_2$.

\noindent
{\bf A.}  $P_2$ contains both $\bot$ and $\top$:
$$T_{P_2}=\ \ \lower 0.4in\hbox{$\epsfysize=0.9in
 \epsfbox{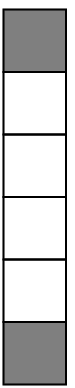}$}\ .
$$

\noindent
{\bf A-I.} $P_1$ has also contains $\bot$ and $\top$:
$$T_P=\ T_{P_1*P_2}=\ \ \lower0.4in\hbox{$\epsfysize=0.9in
\epsfbox{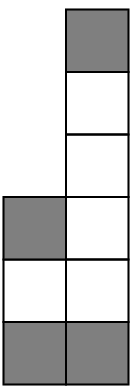}$}$$
\noindent
so $\nu_1-\nu_n=k$ and $\eta_1-\eta_n=k$. In this case define $\Phi(P)=\ \Psi(P)$. Again $\Phi$
is well-defined because the $\bot$-labels in $P_1$ and $P_2$ will actually be consecutive
letters in the word of $P=P_1*P_2$, so they will be paired with each
other. As a result the $\bot$-label of the second column will not move into
the first column. This means that the first row of $\eta'$ has the same
length
as the first row of $\eta$ and therefore
$\eta'\in \Pi^{(n,k)}$.
The image has the form 
$$T_{\Phi(P)}=\ \ \lower0.4in\hbox{$\epsfysize=0.75in
\epsfbox{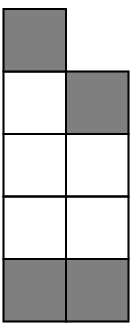}$}$$
\noindent
which is not a column strict tableau.

\noindent
{\bf A-II.} $P_1$ only contains $\bot$:
$$T_P=\ \ \lower0.4in\hbox{$\epsfysize=0.9in
\epsfbox{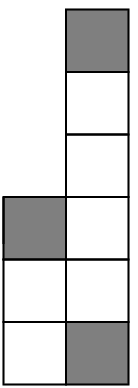}$}$$
\noindent
so $\nu_1-\nu_n=k\,,\eta_1-\eta_n=k$ and $\lambda_1-\lambda_n=k-1$. Using the same argument as
before $\Phi(P)=\ \Psi(P)$ is well-defined. Since $\bot$ in $P_2$ will
not move it follows that
the image
$$T_{\Phi(P)}=\ \ \lower0.4in\hbox{$\epsfysize=0.75in
\epsfbox{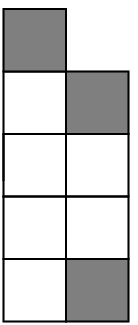}$}$$
\noindent
is not a column strict tableau.

\noindent
{\bf A-III.} The case when $P_1$ only contains $\top$:
$$T_P=\ \ \lower0.4in\hbox{$\epsfysize=0.9in
\epsfbox{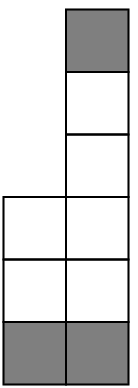}$}$$
\noindent
is not possible as we have $\nu_1-\nu_n=k\,,\eta_1-\eta_n=k$ and $\lambda_1-\lambda_n=k+1$, so $\lambda \notin
\Pi^{(n,k)}$.

\noindent
{\bf A-IV.} $P_1$ does not contain either $\bot$ or $\top$:
\begin{equation}\label{eti}
T_P=\ \ \lower0.4in\hbox{$\epsfysize=0.9in 
\epsfbox{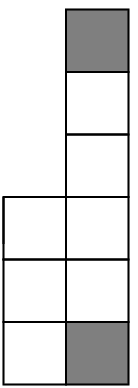}$}
\end{equation}
This is a case when, by applying the operator $\Psi$, it is possible that $\bot$ from the
second column will move into the first column, and hence the possibility
that $ \eta' \notin \Pi^{(n,k)}$.
The operator that we therefore need to construct here will be a  modification of the
operator $\Psi$.
There are two subcases to consider depending on whether $\bot$ is a paired parenthesis in the
word of $P$ or not. 

{\bf a). } If $\bot$ is paired then $\Phi(P)=\Psi(P)$ is well defined since the
pairing of this $\bot$-label means that it remains in the second column i.e.\ the intermediate image partition 
$\eta'$ will have $\eta_1'=\eta_1$, and hence $\eta' \in \Pi^{(n,k)}$. 

{\sl Example 2: If $\nu/\lambda$ = \lower 0.4in\hbox{$\epsfxsize=0.6in
\epsfbox{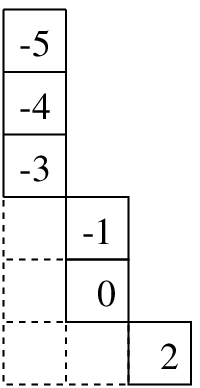}$}, $n=6$, $k=2$, and the tableau $T_P$ =\
\lower0.2in\hbox{$\epsfxsize=0.4in \epsfbox{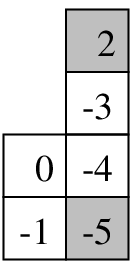}$} , then
\indent $w=\bar{5}\,\bar{4}\,\bar{3}\,\bar{1}\,0\,2$, and its
parentheses  structure is $ \hat{w}=)))(()$. Thus $\widehat{\Phi(w)}=))((()$ and}

$$
T_{\Phi(P)}=\ \lower0.2in\hbox{$\epsfxsize=0.4in
  \epsfbox{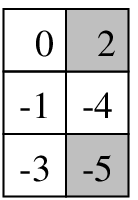}$}\ . 
$$
\noindent
The image has the form
$$T_{\Phi(P)}=\ \ \lower0.4in\hbox{$\epsfysize=0.75in
\epsfbox{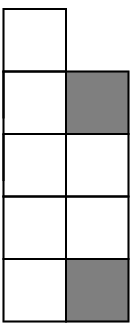}$}$$
\noindent
and  $T_{\Phi(P)}$ is not a column strict tableau.

{\bf b). } The $\bot$-label is not a paired parenthesis in the word of $P$.
This is a case when
$ \Psi$ cannot be applied since $\bot$ would move into the first
column, which means
that $\eta_1'=\eta_1+1$ and $\eta_n'=\eta_n$, so
$\eta_1'-\eta_n'=k+1$. For this case we will define a new operator. 

\noindent
We denote by $\mathcal{D}_1 \subseteq \mathcal A$ the subset of paths
satisfying

\begin{enumerate}

\item[i)] $\nu$ is an edge diagram ($\nu_1-\nu_n=k$)

\item[ii)] the $(\bot,\,\top)$-structure as described by~(\ref{eti})

\item[iii)] $\bot$ is not a paired parenthesis.

\end{enumerate}

\noindent
We write the word of $P$ in the same manner as before and assign parentheses. The first letter in
this word is actually $\top$ from $P_2$ and the last letter is 
$\bot$ from $P_2$, since these numbers are the smallest and the
largest of all labels, respectively i.e.
 $$
w\ =\ \hbox{$\top$............ $\bot$ }.
$$
Suppose that in this word we have $w=\ldots \ldots a_1\ldots a_2\ldots
a_s$, where $a_s=\bot$, and the sequence $a_1a_2\ldots a_s$ represents
the labels corresponding to the last
column of $\nu/\lambda$.

{\sl Example 3:}
$$
\nu/\lambda\ =\lower0.45in\hbox{$\epsfxsize=1in
\epsfbox{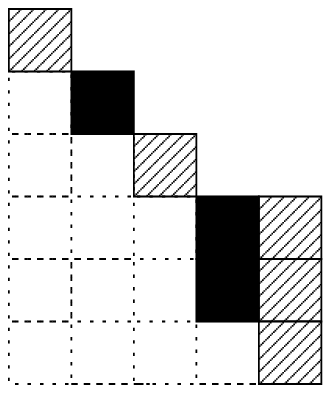}$}\rlap{\hskip 3pt\hbox{$a_1$}}
\rlap{\hskip 3pt\lower0.2in\hbox{$a_2$}}\hskip
3pt\lower0.4in\hbox{$a_3\,\leftarrow\,\bot$}
$$

{\sl where
\lower0.07in\hbox{$\epsfxsize=0.2in\epsfbox{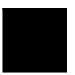}$}
represents boxes in $P_1$ and 
\lower0.07in\hbox{$\epsfxsize=0.2in\epsfbox{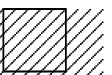}$}
represents boxes in $P_2$.}
\vskip 12pt
Let $i_0= min\{i\ |\ \hbox{$a_i$ is unpaired
letter}\}$.
We note that in this case all left parentheses will be paired, since
the first and last letter in the word $w$ are right unpaired
parentheses. Thus $w$ has the following parentheses structure.
$$
\begin{array}[t]{lcrl}
\hat{w} & = & \ )\,)\,(\,)\,)\,(\,)....(\,(\,)\ \ )\ \
\boldsymbol{\Big )}\ \,) &  \\
w&= \,&                         a_1\,a_2\,a_{i_0}\ & .
\end{array}
$$
\noindent
In the above description of the word $w$, we highlighted the parenthesis
associated to label $a_{i_0}$.
Since all left parentheses (in $P_1$) are paired, the number of
unpaired  (right) parentheses is $|P_2|-|P_1|$.

\noindent
We define the operator
$\phi_1:\ \mathcal{D}_1 \rightarrow\ \mathcal{B}$  on $w$ by
specifying the changes with respect to the parentheses structure so
$\widehat{\phi_1(w)}$=changes all right unpaired parentheses into left
parentheses except the label $a_{i_0}$:
\begin{equation}
\widehat{\phi_1(w)}=\ ((()(()\ldots (())\boldsymbol{\Big )}( .
\end{equation}\label{newinv}

{\sl Example 4:}

{\sl For $\nu/\lambda$=\ 
\ \lower0.4in\hbox{$\epsfxsize=0.8in\epsfbox{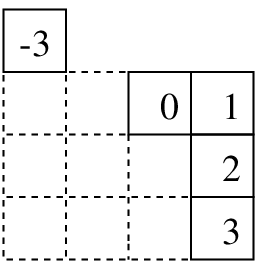}$},
$n=4$, $k=3$ and
$T_P$=
\lower0.3in\hbox{$\epsfxsize=0.4in\epsfbox{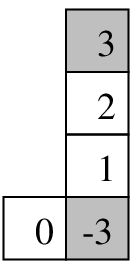}$} we have
$w=\bar{3}\,0\,1\,2\,3$
and the \indent parentheses structure is $\hat{w}=)()\boldsymbol{\Big )})$. Thus
$\widehat{\phi_1(w)}=(()\boldsymbol{\Big )}($ and $T_{\phi_1(P)}$=\ 
\lower0.2in\hbox{$\epsfxsize=0.4in\epsfbox{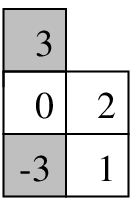}$} .}
\vskip 12pt
\noindent

In the case when $a_{i_0}=a_s=\,\bot$ the situation is slightly different.

{\sl Example 5:}

{\sl For $\nu/\lambda$=\ 
\ \lower0.2in\hbox{$\epsfxsize=0.6in\epsfbox{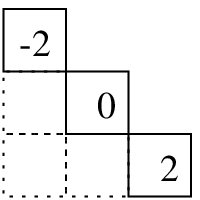}$},
$n=3$, $k=2$ and
$T_P$=\ \ 
\lower0.15in\hbox{$\epsfxsize=0.4in\epsfbox{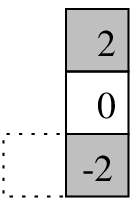}$} we have
$w=\bar{2}\,0\,2$
and the parentheses \indent structure is $\hat w=))\boldsymbol{\Big )}$. Thus
$\widehat{\phi_1(w)}=((\boldsymbol{\Big )}$ and $T_{\phi_1(P)}$=\ 
\lower0.15in\hbox{$\epsfxsize=0.4in\epsfbox{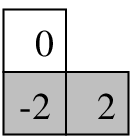}$} .}
\vskip 12pt
 As we have seen in {\sl Example 4 and 5} the image has the form
\begin{equation}\label{image}
T_{\phi_1(P)}=\ \
\lower0.4in\hbox{$\epsfysize=0.75in
\epsfbox{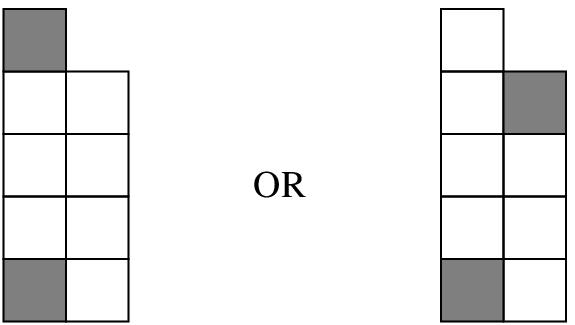}$}\ .
\end{equation}

\begin{prop}
The operator $\phi_1$ is well defined.
\end{prop}
\noindent
{\bf Proof.} \

There are two things that we need to check. Given that
$\phi_1(P)=P_1'*P_2'$ with intermediate diagram $\eta'$, so $\lambda\ \buildrel P_1' \over \longrightarrow\ \eta' \
\buildrel P_2' \over \longrightarrow\ \nu$, we have to see that
\begin{itemize}
\item $\eta'$ is a partition

\item the skew shapes $\eta'/\lambda$ and $\nu/\eta'$ are
  column-strips.
\end{itemize}
 
We show first that $\eta'$ is a partition.
Assume that $\eta'$ is not, that is there exists $l$ such
that $\eta_l'<\eta_{l+1}'$. Since the operator $\phi_1$ removes labels from a column strip we must have
$$
 \eta_{l+1}'=\eta_l'+1
$$
$$
 \epsfbox{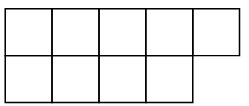}\ .
$$
For simplicity let us denote the labels in the first column of $P$ and
$\phi_1(P)$ by 
$1$ and the ones in the second column by $2$. Generically $\phi_1$
transforms some ``$2 \rightarrow 1$''.

\noindent
We obtain the above situation only if
$$
\lower 0.1in \hbox{$\epsfbox{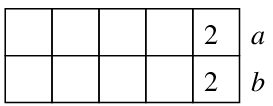}$}\ \buildrel \phi_1
\over \longrightarrow\ \ \lower 0.1in
\hbox{$\epsfbox{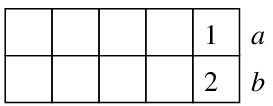}$}\ .
$$
\noindent
This means that the label $a$ is not a paired parenthesis in $w$ and the label
$b$, which is a right parenthesis, is paired or is the label
$a_{i_0}$. Let us consider these two situations.

--If $b\,(\sim 2)$ is paired, then its pair, a label $b'\,(\sim 1)$, must be
to its left in $w$. Two situations may occur

\noindent
1).  $b'\leq a$. We have
$$
\begin{array}[t]{lcl}
    \hat{w}=...(... & ) & )... \\
    w=...b' & a & b...\ .
\end{array}
$$ 
\noindent
which shows that the label $a$ would be paired. This is not
possible since $a$ is assumed to be unpaired.

\noindent
2). $b'>b$. Since $b'$ must be on $b$'s left we have $b'=b$.
$$
\begin{array}[t]{lcl}
    \hat{w}=...) & ( & )... \\
         w=...a & b' & b...\ .
\end{array}       
$$
\noindent
In the shape $\nu/\lambda$  this corresponds to
$$
\epsfbox{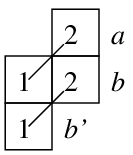}\ .
$$
\noindent
In the figure we also indicated the pairing. The pairing that we
illustrated above is a consequence of the fact that labels on the same
diagonal are in fact equal so they are consecutive letters in the word.
In this case we note that there must exist a label from the first path
($\sim 1$) that is above $b'$. But since this label and $a$ are on the
same diagonal (equal) they will pair, a
contradiction.

--Therefore assume $b$ is $a_{i_0}$ (the special label). In this case the labels $a$
and $b$ are part of the last column of $\nu$, and  $a=a_{i_0}-1,\,b=a_{i_0}$.
Since $a_{i_0}$ was defined to be the smallest label in the last
column to be unpaired, the smaller labels in this column $a_i<a_{i_0}$
are paired parentheses. Thus, in particular $a$ would be paired and we again
obtain a contradiction. Hence $\eta'$ is indeed a partition. 
\vskip 12pt
We now prove that $\eta'/\lambda$ and $\nu/\eta'$ are column
strips.

\noindent
Since $\phi_1$ moves labels from the second path into the first path we
have
$$
 \nu/\eta' \subset \nu/\eta
$$
and because  $\nu/\eta$ was a column strip, $\nu/\eta'$ is a column
strip  as well.

Next we show that $\eta'/\lambda$ is a column strip.
Assume it is not. This occurs when
in two consecutive columns and the same row in $T_P$ we have first
a label $a\,(\sim 1)$ followed by $b\,(\sim 2)$ changed by 
$\phi_1$ into two $1$'s i.e.
$$
\lower 0.1in \hbox{$\epsfbox{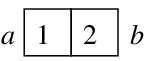}$}\ \buildrel \phi_1
\over \longrightarrow\ \ \lower 0.1in
\hbox{$\epsfbox{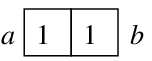}$}\ .
$$

\noindent
Let us study the situation in $T_P$. We note that in $\hat{w}$ the label $a$ is paired (since it is a left
parenthesis and all of them are paired) but the label $b$ is not a
paired parenthesis and it is not $a_{i_0}$.
The situation above has the following features.

--There is no other label $1$ (in the first path) underneath the label
$a$. If there were any, say
$$
\epsfbox{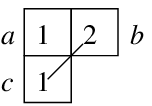}
$$
\noindent
the label $c\,(\sim 1)$ just below $a$ would pair with $b$.

--There must be other labels from the second path, ($2$'s) above the
label $b$. If there were none, then $a$ would be paired with
$b$ in $\hat{w}$ since $b$ is the first right parenthesis on the right of
$a$ i.e.
$$ 
\begin{array}[t]{lrll}
      \hat{w}=&....( & )... &  \\
           w=& a & b\ & .
\end{array}
$$
We observe that the number of $1$'s in the first column above the
label $a$ exceeds or is equal to the number of $2$'s in the
second column above $b$ i.e.
$$
\epsfbox{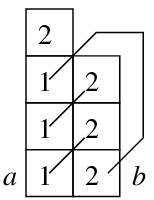}\ .
$$
\noindent
Now since the labels on diagonals are the first to pair, the label $b$ will
pair with the label $1$ in the first column situated on the same row as
the last $2$ in the second column. This again contradicts our
requirement that $b$ is not paired.
 Thus $\eta'/\lambda$ is a column strip.

\noindent
This concludes the proof of Proposition 8.
\hfill
$\Box$

 Our next task is to find a complete characterization of the image
 inside $\mathcal B$ and
 to define $\Phi$ in this case. 
 Consider $w$ to be the word of $T_P$, for $P\in \mathcal B$  and let $w'= a_1a_2\ldots a_s$ be the subword of $w$ made with the labels in
the last column of the partition $\nu$ so 

 $$
 \begin{array}[t]{lcrl}

\hat{w} & = & (\,(...(\, (\, )\,(\,) \,)\,)\, )\, (\, (\, (\ \ &   \\
\widehat{w'} & = &         )\,\,\,)\,)\,)\,)\, (\, (\, (\ \  &  \\
w' & = &         a_1\, a_2.......a_s & .
 \end{array}
 $$
This subword $w'$ might contain labels from both $P_1$ and
$P_2$. Since this is a word of a column we must first have the 
boxes from $P_1$  on top of which are the boxes from $P_2$.

{\sl Example 6:} Consider
$$
\nu\ =\lower0.45in\hbox{$\epsfxsize=0.65in
\epsfbox{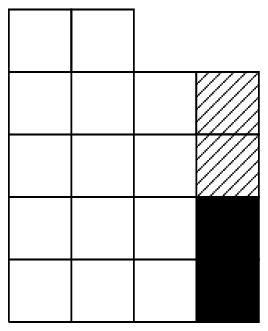}$}\rlap{\hskip 3pt\hbox{$\Big{\}}$
in $P_2$}}
\hskip 3pt\lower0.3in\hbox{$\Big{\}}$ in $P_1$}
$$
\noindent
Thus the parentheses structure is 
$$
\begin{array}[t]{rl}
\widehat{w'}= & )\, )....)\,\boldsymbol{\Big )}(\, (.... (\\
   w= &  a_1..... a_{i_0} 
\end{array}
$$
\noindent i.e. the right are
followed
by the left parentheses.  We denote by $a_{i_0}$ the rightmost right
parenthesis of $w'$. We also identify the
label $(a_1-1)$ (if it exists in $w$) which will play a role in the next
definition. This is the label situated on the penultimate column and
on the same row as the last box in the last column, i.e. the label $a_1$.

{\sl Example 7:}

$$
\nu/\lambda\ =\lower0.45in\hbox{$\epsfxsize=1in
\epsfbox{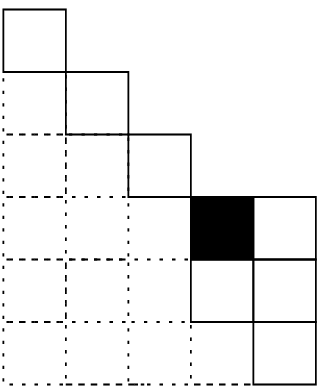}$}
\rlap{\hskip -0.4in\lower-0.2in\hbox{($a_1-1$)}}
\rlap{\hskip 3pt\hbox{$a_1$}}
\rlap{\hskip 3pt\lower0.2in\hbox{$a_2$}}\hskip
3pt\lower0.4in\hbox{$a_3\,\leftarrow \bot$}
$$
\centerline{{\bf Figure 3.} {\sl The label $(a_1-1)$}.}
\vskip 12pt
{\bf Definition.} 
The subset $\mathcal{D}_2 \subseteq \mathcal{B}$ of paths $P=P_1*P_2$
is defined by the paths satisfying
\begin{enumerate}\label{set2}

\item $T_P$ is a column strict tableaux

\item the ($\bot,\, \top$)-label structure of $T_P$ is described
  by~(\ref{image})

\item the last column contains labels from $P_2$ and the label $a_{i_0}$ is not paired with $(a_1-1)$ (if the
  latter exists in $w$)

\item the smallest label $\top$ in $w$, is either an
unpaired left parenthesis or paired with the label $a_{i_0}$.

\end{enumerate}

\noindent
The operator $\phi_2:\ \mathcal{D}_2\rightarrow
\mathcal{A}$ applied to $w$ is defined via

$\widehat{\phi_2(w)}$ which changes all unpaired left parentheses
into right ones including the left parenthesis $b_{i_0}$ paired with
$a_{i_0}$ i.e.
$$
 \begin{array}[t]{lcrl}
w & = & .......b_{i_0}......... .a_{i_0}... &             \\
\hat{w} & = & (\,(\,(\,)\,(\,(\, (\,)\,(\,) \,)\,\Big{)}\, (\, (\, &  \\
\widehat{\phi_2(w)}& =
&)\,)\,(\,)\,)\,(\,(\,)\,(\,)\,)\,\Big{)}\,)\,)\, & .     

 \end{array}
 $$

\noindent
\begin{rem}
The following properties of the image $T_{\phi_2(P)}$ are easy
consequences of the above definition.

\begin{enumerate}

\item[i)] The number of parenthesis to be changed is now
  $|P_1|-|P_2|-1$. If $|P_1|=p$, $|P_2|=q$ then
  $T_{\phi_2(P)}=P_1''*P_2''$ with the intermediate diagram $\eta''$:
$\lambda\ \buildrel P_1'' \over \longrightarrow\ \eta'' \
\buildrel P_2''\over \longrightarrow\ \nu$ and $|P_1''|=p-(p-q+1)=q-1$, $|P_2|=q+p-q+1=p+1$.

\item[ii)] The first letter $\top$ in $w$, is moved by  the
  operator $\phi_2$ in $P_2''$.
     
\item[iii)] The biggest letter $\bot$ in $w$, is either an
  unpaired label or is the label $a_{i_0}$. In both cases this label is in $P_2''$ and is unpaired in $\phi_2(w)$.

\end{enumerate}

\noindent
These characteristics prove that $\hbox{Im}(\phi_2) \subseteq \hbox{Dom}(\phi_1)=\mathcal{D}_1$.

\end{rem}
\begin{prop}
The operator $\phi_2$ is well defined.
\end{prop}

\noindent
{\bf Proof.}

\noindent
As in Proposition 8 we have to check that the intermediate diagram
$\eta''$ is a partition and that both skew diagrams
$\eta''/\lambda$ and $\nu/\eta''$ are column strips.
The proof that $\eta''$ is a partition is similar to the one in
Proposition 8 and we leave it to the reader.

We shall prove that $\eta''/\lambda$ and $ \nu/\eta''$ are
column-strips.
Since $\phi_2$ moves labels from the first path into the second path we
have that $\eta''/\lambda \subset \eta/\lambda$ so $\eta''/\lambda$ is a column strip.

We now show that $\nu/\eta''$ is a column strip. Assume it is not. This
occurs when in two consecutive columns and the same row in $T_P$ we first
have a label $a(\sim 1)$ followed by $b(\sim 2)$ changed by $\phi_2$
into two $2$'s i.e.
$$
\lower 0.1in \hbox{$\epsfbox{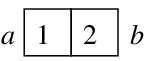}$}\ \buildrel \phi_2
\over \longrightarrow\ \ \lower 0.1in
\hbox{$\epsfbox{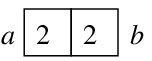}$}\ .
$$
\noindent
This means the label $a$ is either an unpaired left parenthesis or is
the label $b_{i_0}$.

--Assume that $a$ is an unpaired left parenthesis.
We first note that there are no other $2$'s above the label $b$ (if there
  were any, the first label $2$ above $b$ would be on the same diagonal
with $a$, so it would pair with it).
Another useful observation is that the number of $2$'s in the second
column must exceed or be equal to the number of $1$'s in the first
column below the label $a$ i.e. 
$$
\epsfbox{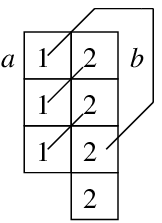}\ .
$$
\noindent
In this case the label $a$ pairs with the label $2$ in the second
column situated on the same row as the first $1$ in the first column.

--Assume that $a$ is the label $b_{i_0}$, i.e. the label paired with
$a_{i_0}$. As before we claim that there are no labels $2$ in the
second column above $b$. If there were any, the first one above $b$
would pair
with $a$, so this label must be $a_{i_0}$. This is not possible since
there are no labels $2$ below $a_{i_0}$, by the definition of $a_{i_0}$. As above we also have that the number of $2$'s in the
second
column must exceed or be equal to the number of $1$'s in the first
column below the label $a$. We note that $a=b_{i_0}$ pairs with the label
$a_{i_0}$ ($\sim 2$) in the second column situated on the same row as the
first label $1$ in the first column.
    This shows that the value of the label $a=b_{i_0}$ is $a_1-1$,
    where $a_1$ is the last label in the last column of
    $\nu$, which is also the label $b$.

\noindent
This situation contradicts condition (3) in the definition of 
$\mathcal{D}_2$.

\noindent
This concludes the proof of Proposition 10.
\hfill
$\Box$
\vskip 12pt

The following result shows that the operators $\phi_1$ and
$\phi_2$ are inverse to each other.

\begin{prop}\hfill

a). Im($\phi_1)\subset \mathcal{D}_2$ and Im($\phi_2) \subset
\mathcal{D}_1$.

b). $\phi_1 \circ \phi_2= id_{\mathcal{D}_2}$ and $\phi_2 \circ
\phi_1= id_{\mathcal{D}_1}$.

\end{prop}

\noindent
{\bf Proof.}

a). In Remark 9 we showed that Im($\phi_2) \subset \mathcal{D}_1$. We next
show that Im($\phi_1)\subset \mathcal{D}_2$.

\begin{enumerate}

\item In $\phi_1(w)$ all right parentheses will be paired (including
$a_{i_0}$) so $T_{\phi_1(P)}$ is a column strict tableau.

\item We also establish in description~(\ref{image}) the ($\top\,,\bot$)-label structure.

\item The label $a_{i_0}$ in $\phi_1(w)$ cannot pair with the label
$(a_1-1)$ (see {\bf Figure 3}).
If this happens, then in $w$ the label $(a_1-1)$ was a right
parenthesis i.e in the second path. However in $w$ the label $a_1$, which
is situated on the same row, is also in the second path. This cannot
be possible since the second path must represent a column-strip.

\item The $\top$-label, which is the first letter in $w$ (or $\phi_1(w)$)
is unpaired or it pairs with $a_{i_0}$ if there were no other unpaired right parentheses
between $\top$ and $a_{i_0}$.

\end{enumerate}

\noindent
Hence Im($\phi_1)\subset \mathcal{D}_2$.

b). The relation $\phi_2 \circ \phi_1= id_{\mathcal{D}_1}$ is obvious
by the definition of the operators $\phi_1$ and $\phi_2$. We
illustrate this by the following example. Let
$T_P\in \mathcal{D}_1$ and $w$ be its word with the 
parentheses structure
$$
\hat{w}\ =\ )))()(())))(()())\Big{)}))
$$
\noindent
where $a_{i_0}$ is highlighted. Applying $\phi_1$ we get.
$$
\widehat{\phi_1(w)}\ =\ (((()(())(((()())\Big{)}((.
$$

\noindent
Since we showed that $T_{\phi_1(P)}\in \mathcal{D}_2$ we can apply the
operator $\phi_2$ to $\phi_1(w)$ to get
$$
 \widehat{\phi_2(\phi_1(w))}\ =\ )))()(())))(()())\Big{)})).  
$$
\noindent
Therefore we have $\phi_2(\phi_1(w))=w$.
Similarly we have that $\phi_1 \circ \phi_2= id_{\mathcal{D}_2}$.

\hfill
$\Box$

\noindent

We define all column strict tableaux of shape $\mu$ that do not belong
to the set $\mathcal{D}_2$ to be {\em $k$-fusion} and we denote their set by
$CSF_k(\mu)= CS(\mu) \setminus \mathcal {D}_2$.

\noindent
We finish the case that we studied (A-IV,b) by letting
$\Phi=\phi_1$.

\noindent
{\bf B.}  $P_2$ contains the $\bot$-label but not the $\top$-label: 
$$T_{P_2}=\ \ \lower0.4in\hbox{$\epsfysize=0.9in
\epsfbox{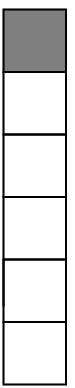}$}\ .
$$

\noindent
{\bf B-I.} $P_1$ contains both labels:
$$T_P=\ \ \lower0.4in\hbox{$\epsfysize=0.9in
\epsfbox{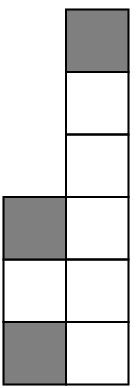}$}\ .
$$
\noindent
In this case $\eta_1-\eta_n=k-1$ and $\lambda_1-\lambda_n=k-1$ so
define $\Phi(P)=\ \Psi(P)$ which is again well defined by a similar
argument to the one in {\bf A}--I. Therefore the tableau of
the image:
$$T_{\Phi(P)}=\ \ \lower0.4in\hbox{$\epsfysize=0.75in
\epsfbox{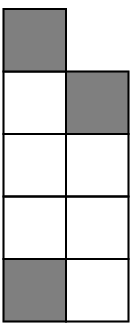}$}$$
\noindent
is not a column strict tableau.

\noindent
{\bf B-II.} $P_1$ contains the $\bot$-label but not the $\top$-label: 
$$T_P=\ \ \lower0.4in\hbox{$\epsfysize=0.9in
\epsfbox{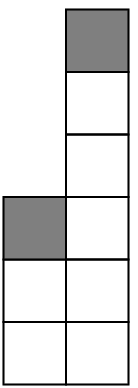}$}\ .$$
\noindent
Here it is clear that $\Phi(P)=\ \Psi(P)$ is well-defined and the
image 
$$T_{\Phi(P)}=\ \ \hbox{$\epsfysize=0.75in
\epsfbox{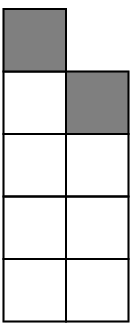}$}$$
\noindent
is again not a column strict tableau.

\noindent
{\bf B-III.} $P_1$ contains the $\top$-label but not the $\bot$-label:
$$T_P=\ \ \lower0.4in\hbox{$ \epsfysize=0.9in
\epsfbox{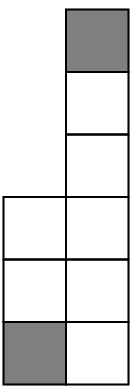}$}\ .$$
\noindent
In this case $\eta_1-\eta_n=k-1$ and $\lambda_1-\lambda_n=k$. Define $\Phi(P)=\ \Psi(P)$. Here it
is possible that the $\bot$-label will move into the first column if it is not a paired
parenthesis, which means that the length of the first row of $\eta '$
will increase by one, so
$\eta_1 '-\eta_n '=(k-1)+1=k$, and $\eta'$ is still a
$(n,k)$-restricted partition. In this case the image can be 
$$T_{\Phi(P)}=\ \ \lower0.4in\hbox{$\epsfysize=0.75in
\epsfbox{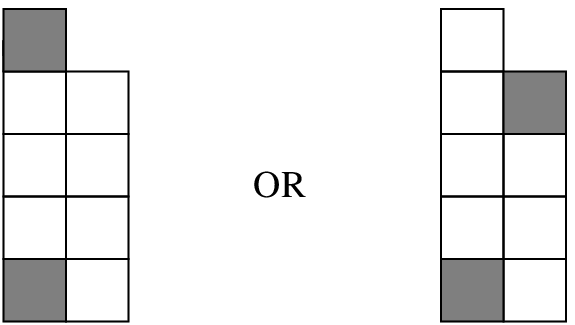}$}\ .$$
\noindent
In both cases the image is not a column strict tableau.

\noindent
{\bf B-IV.} $P_1$ does not contain either labels:
$$T_P=\ \ \lower0.4in\hbox{$\epsfysize=0.9in
\epsfbox{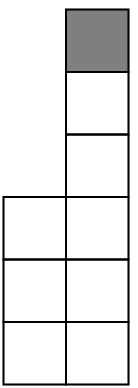}$}\ .$$
\noindent
This is a case similar to the previous one, so $\Phi$ will be defined in the same way. The
image can be
$$T_{\Phi(P)}=\ \ \lower0.4in\hbox{$\epsfysize=0.75in
\epsfbox{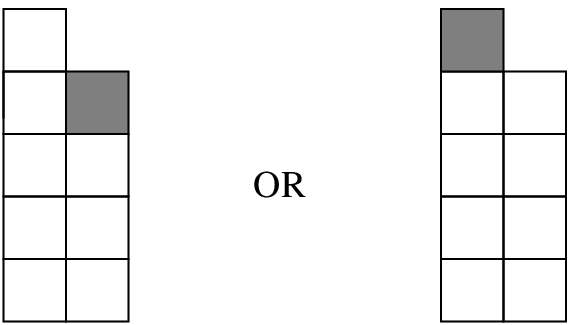}$}$$
\noindent
but, again, in both cases the image is not a column strict tableau.

For the remaining cases {\bf C} where $P_2$ has the $\top$-label but not the $\bot$-label and {\bf D}
when $P_2$ does not have either labels, since the $\bot$-label is not present there is no danger
in increasing the first row of the diagram $\eta '$, so in all these
cases we define $\Phi(P)=\ \Psi(P)$. The structure of the ($\bot,\, \top$)-label of the images will look the same as
for $P$ and the associated tableaux are not column strict.  

{\sl Observation:} As we have seen in {\bf A}--I, many of these 8 remaining
cases will not be possible.
To conclude we have found only one case where we introduce a new operator.

To finish this analysis we must also consider the situation
$P=P_1*P_2$ for which $|P_1|\geq |P_2|$.
Since we have already seen the structure of the image of $\Phi$, we have the following.

{\bf i)} If there is a violation of the $CS$-property for $T_P$ and
$P \notin \mathcal {D}_1$ then $\Phi(P)=\Psi(P)$. 

{\bf ii)} If $P \in \mathcal{D}_1$, then $\Phi(P)=\phi_1(P)$.

{\bf iii)} If $T_P \in CS(\mu)\setminus CSF_k(\mu)$, i.e. $P \in
\mathcal{D}_2$ then $\Phi(P)=\phi_2(P)$.

{\bf iv)} In any other case, i.e. $T_P \in CSF_k(\mu)$, we have $\Phi(P)=P$. 

\noindent
In fact we also proved that $\Phi$ is an involution on the set of paths made from two decreasing
paths whose intermediate partitions are $(n,k)$-restricted.

The fusion coefficients, which are the number of fixed points of the
involution $\Phi$, count the number of {\sl $k$-fusion} tableaux. Therefore we have:
\begin{thm}
For $\mu$ a two-column partition and any level $k$ we have
\begin{equation}\label{eq:fus}
  N_{\lambda \mu}^{(k)\nu}= \sharp \{\hbox{paths $P$ from $\lambda$ to
$\nu$ in $\mathcal{P}_{(\mu')}^{(n,k)}$ that fit $\mu$ and 
$T_P \in CSF_k(\mu)$}\}.
\end{equation}
\end{thm} 

\begin{rem}
 By replacing every label from $P_1$ by $1$ and the labels from $P_2$
 by $2$ in the partition $\nu$ and reinterpreting the conditions in
 the definition of $\mathcal{D}_2$ we get the following characterization
 for fusion coefficients in the case $\mu$ has two columns.

 The coefficient $N_{\lambda \mu}^{(k)\nu}$ counts the number tableaux in
    $RS\Pi^{(n,k)}(\nu/\lambda,\, \mu)$ whose word (read
    column-wise) is lattice, except the tableaux for which
\begin{itemize}  

\item  $\nu_1-\nu_n=k$,

\item  the first row contains exactly one of $1$ or $2$ and the last row
contains exactly a $1$, 

\item the last column contains $2$'s,

\item the number of $1$'s in the penultimate column under the height
  of the last column is strictly less than the number of $2$'s in the
  last column (see {\bf Figure 4} under the thick line),

\item  the number of $1$'s in the reading
  word is always strictly bigger than the number of $2$'s
  except (perhaps) when the last $2$ is counted.

\end{itemize}

$$
\epsfbox{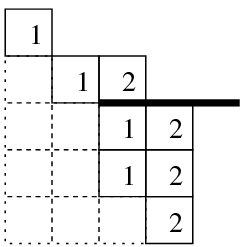}
$$
\centerline{{\bf Figure 4. }{\sl Example of a tableau described above
    for $n=5$ and $k=3$.}}
\end{rem}

\section{ Applications}

We shall now give some consequences of the last theorem.
\begin{cor}[part of Prop(2.2) of~\cite{GdWz}]
 For any level $k$, if $\mu$ is a one or two-column partition, we have

a). $N_{\lambda \mu}^{(k)
\nu}\leq c_{\lambda
\mu}^{\nu}$, where the latter are the classical Littlewood-Richardson coefficients.

b). If all the paths in $\Omega_k$ are only
passing through $(n,k)$-restricted partitions (e.g.
$\lambda_1-\lambda_n\leq k-1$), then $N_{\lambda \mu}^{(k) \nu}=c_{\lambda
\mu}^\nu$.

\end{cor}
\noindent
{\bf Proof.}

a). This is an obvious consequence of Theorem 12.

b). With this condition, the case where $T_P \notin CSF_k(\mu)$ cannot
occur, so $\Phi=\Psi$.
\hfill
$\Box$
\vskip 12pt  
The next result proves the conjecture~(2.4) in~\cite{W2}, in our special case.

\begin{thm}
If $\mu$ is a one or two-column partition, then we have
$$
N_{\lambda \mu}^{(k) \nu} \leq N_{\lambda \mu}^{(k+1) \nu}.
$$
\end{thm}

\noindent
{\bf Proof.}

From the way we constructed the involution $\Phi$ we know that 
$$
N_{\lambda \mu}^{(k) \nu}=\sharp \{P \in \mathcal{P}_{(\mu')}^{(n,k)},\ \hbox{ from $\lambda$ to $\nu$  such that $\Phi(P)=P$}\}.$$
\noindent
In order to prove the inequality it suffices to see that
if
$\Phi_{(k)}(P)=P$ then
$\Phi_{(k+1)}(P)=P$ as well. Note we index the operators by the levels that we
consider.

\noindent
Suppose this is not the case. Since $T_P \in CSF_k(\mu)\subset CS(\mu)$, it is
possible that $P$ is not $(k+1)$-fusion. In this case the
partition $\nu$ must be an edge
diagram for level $(k+1)$ i.e. $\nu_1-\nu_n=k+1$. This cannot happen since
$\nu$
is also a $(n,k)$-restricted partition i.e. $\nu_1-\nu_n \leq k$.
\hfill
$\Box$ 
\vskip 12pt
We yield another application by using the rank-level duality.
Recall~\cite{GdWz} that we can define a bijection between
$(n,k)$-restricted partitions and $(k,n)$-restricted partitions as
follows.

\noindent
For $\lambda \in \Pi^{(n,k)}$, cut the rectangle $\lambda_1 \times n$
into rectangles of sides $k \times n$. Conjugate each rectangle separately and then glue the resulting
partitions back together.

{\sl Example:}
$$
\lower0.3in\hbox{$\epsfxsize=1.1in \epsfysize=0.9in
\epsfbox{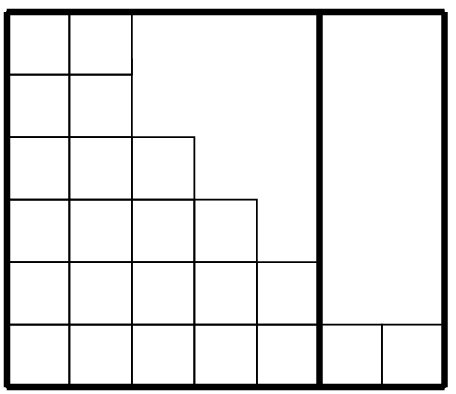}$}\ \ \ \ 
\boldmath{\rightarrow}\ \ \ \
\lower0.3in\hbox{$\epsfxsize=1.1in \epsfysize=0.8in
\epsfbox{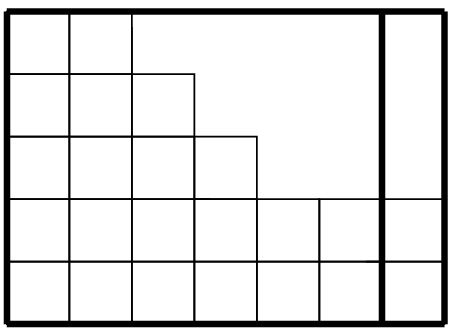}$}\ .$$
\noindent
It is clear that the resulting partition $\tilde{\lambda}$ constructed in
this way is a $(k,n)$-restricted partition. Goodman and Wenzl~\cite{GdWz}
showed that
the fusion coefficients are invariant under this bijection i.e.
$N_{\lambda \mu}^{(k)\nu}=N_{\tilde{\lambda}
\tilde{\mu}}^{(n)\tilde{\nu}}$ and as a result we have the following theorem.
\begin{thm}
For any level $k$ if $n\geq 3$, and $\mu$ is a partition with one or
two rows, then the
fusion
coefficients $N_{\lambda\mu}^{(k)\nu}=N_{\tilde{\lambda}\tilde{\mu}}^{(n)\tilde{\nu}}$ are given
in Theorem 12.
\end{thm}
\noindent
{\bf Proof.}
If $n \geq 3$ then $\mu_1\leq k$. Therefore $\tilde{\mu}=\mu'$, where
$\mu'$ is the conjugate of $\mu$. It is now clear that we are in the
setting of Theorem 12 and as a result we can determine the
coefficients $N_{\lambda\mu}^{(k)\nu}$.
 
\hfill
$\Box$
\begin{rem}
For $n=2$ the fusion coefficients are given by the Gepner-Witten
formula~\cite{GepWi}
\begin{equation}\label{GeWi}
N_{\lambda \mu}^{(k)\nu}=\left\{ \begin{array}{ll}
                        c_{\lambda \mu}^{\nu} & \hbox{ if } k\geq
                        \lambda_1-\lambda_2+\mu_1-\mu_2+\nu_1-\nu_2\\
                        0 & \hbox{ otherwise }.
                        \end{array}
                        \right.
\end{equation}
\end{rem}

\section {Conclusions}

The goal of this paper was to find an appropriate involution, in the same 
manner as for the Littlewood-Richardson rule, which would give a much
desired combinatorial description for the fusion coefficients.
As for the LR-rule we started by defining the involution in the case 
where one partition has at most two columns. We were able to prove that
except in one case, the involution remained the same. In this special
case we argued that we must construct a different operator, somehow
similar with the classical one, and we were successful in doing so.
The obstruction in defining the involution in the general case is the fact that we could not find a canonical position in the partition
where the operator for the $2$-column case is to be applied. A reason for
this is that it seems there is no specific area in the $2$-column part
that remained unchanged by the involution. 

Another question one can ask about fusion coefficients is does there
exist an equivalent Robinson-Schensted
correspondence? This question seems legitimate
since we can establish a result similar to the following equality~\cite{BS}.
\begin{prop}
If $\lambda \subseteq \nu$ are two partitions then
\begin{equation}\label{nan}
\sharp\{\hbox{ paths from }\lambda\  \hbox{to}\
\nu\}=\sum_{\mu \vdash |\nu/\lambda|} c_{\lambda \mu}^{\nu}f^\lambda
\end{equation}
where $f^\lambda$ denotes the number of standard tableaux of shape
$\lambda$.
\end{prop}

We extend some of the definitions given previously. We say that a
path
is {\sl $(n,k)$-restricted} if it only passes
through $(n,k)$-restricted partitions. 
Let $T_\lambda$ be a standard tableau of shape $\lambda$. We can
identify the standard tableau with a path from $\O$ to $\lambda$.
The $i^{th}$ partition in the path is obtained obtained from
the previous one by adding the box indexed $i$ in the tableau $T_\lambda$.

{\sl Example: The tableau $T_\lambda$= \ \lower0.15in\hbox{$\epsfxsize=0.5in
  \epsfbox{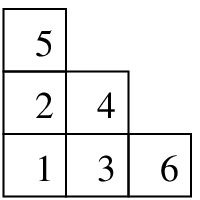}$} represents the path}
$$
\lower0.1in\hbox{$\epsfxsize=0.15in\epsfbox{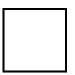}$}
\ \hookrightarrow 
\ \lower0.1in\hbox{$\epsfxsize=0.15in\epsfbox{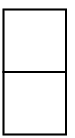}$}
\ \hookrightarrow
\ \lower0.1in\hbox{$\epsfxsize=0.3in\epsfbox{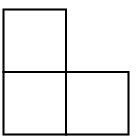}$}
\ \hookrightarrow 
\ \lower0.1in \hbox{$\epsfxsize=0.3in\epsfbox{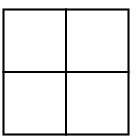}$}
\ \hookrightarrow
\ \lower0.1in\hbox{$\epsfxsize=0.3in\epsfbox{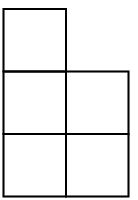}$}
\ \hookrightarrow 
\ \lower0.1in \hbox{$\epsfxsize=0.45in\epsfbox{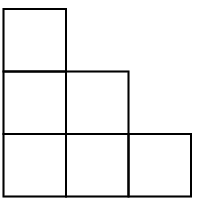}$}\ .
$$
We say that a standard tableau is {\sl $(n,k)$-restricted} if the
associated path is $(n,k)$-restricted. We note that this definition is
consistent with the definition of a column-strict restricted tableau.
We denote by $f_k^\lambda$ the number of
$(n,k)$-restricted standard tableaux of shape $\lambda$. A future task
is to find some sort of expression for these numbers. We have the following result.

\begin{thm} For $\lambda \subseteq \nu$, partitions in $\Pi^{(n,k)}$ we
have
\begin{equation}\label{nice} 
\sharp\{\hbox{restricted paths from }\lambda\  \hbox{to}\
\nu\}=\sum_{\stackrel{\mbox{\scriptsize $\mu \vdash |\nu/\lambda|$}}%
{\mu\in \Pi^{(n,k)}}}  N_{\lambda \mu}^{(k)\nu}f_k^\mu.
\end{equation}
\end{thm}

\noindent
{\bf Proof.}

First we prove by induction that 
\begin{equation}\label{e}
(e_1)^m=\,\sum_{\mu \vdash m}s_\mu f_k^\mu.
\end{equation}
We stress here that all equalities take place in the fusion algebra
$\mathcal{F}^{(n,k)}$ and all partitions involved are $(n,k)$-restricted.

\noindent
If $m=1$ the right hand side of Equation~(\ref{e}) is $s_{(1)}f_k^{(1)}=e_1\cdot 1$.

Assume that the equality is true for $m$ and we shall prove it for
$m+1$. First observe 
$$
(e_1)^{m+1}=(\sum_{\mu \vdash m}s_\mu f_k^\mu)\cdot e_1=\sum_{\mu \vdash
  m} (s_\mu e_1)f_k^\mu= 
\sum_{\mu \vdash m}(\sum_{|\nu/\lambda|=1} s_\nu) f_k^\mu = 
 \sum_{\nu \vdash m+1} s_\nu (\sum_{|\nu/\mu|=1}f_k^\mu).
$$
To show that $\sum_{|\nu/\mu|=1}f_k^\mu =f_k^\nu$ we note that each $(n,k)$-restricted standard tableau
of shape $\mu$ determines a unique $(n,k)$-restricted standard tableau of
shape $\nu$ by adding the corresponding box with the entry
$(m+1)$.
This process is reversible since the box filled with $(m+1)$, which is the
largest number of the standard tableau, is an exterior corner of the shape
$\nu$.  

\noindent
Therefore we get that $(e_1)^{m+1}=\sum_{}s_\nu f_k^\nu$.

Now we proceed to prove Equation~(\ref{nice}). If we multiply 
$s_\lambda$ successively
with $e_1$ in the fusion algebra
we get
\begin{equation}\label{f}
s_\lambda (e_1)^m=\sum_{|\nu/\lambda|=m}\sharp\{\hbox{restricted paths
from }\lambda\  \hbox{to}\ \nu\}s_\nu.
\end{equation}
Using~(\ref{e}), the left-hand side of this equality becomes
$$
s_\lambda (e_1)^m =s_\lambda \sum_{\mu \vdash m} s_\mu f_k^\mu=\sum_{\mu
\vdash m} s_\lambda s_\mu f_k^\mu=\sum_{\mu \vdash m} (\sum_{|\nu/\lambda|=m}
 N_{\lambda \mu}^{(k)\nu}s_\nu) f_k^\mu=\sum_\nu (\sum_{\mu \vdash m}N_{\lambda \mu}^{(k)\nu}f_k^\mu)s_\nu.
$$
Equating the coefficient of $s_\nu$ in the right-hand side of
Equation~(\ref{f}) and the line above we get 
$$     
\sharp\{\hbox{restricted paths from }\lambda\  \hbox{to}\
\nu\}=\sum_{\stackrel{\mbox{\scriptsize $\mu \vdash |\nu/\lambda|$}}%
{\mu\in \Pi^{(n,k)}}}  N_{\lambda \mu}^{(k)\nu}f_k^\mu.
$$
\hfill
$\Box$

In view of the last equation one could hope to define 
fusion-Knuth relations among the words of the restricted paths. This might
happen since in the classical case, the Knuth relations and
Equation~(\ref{nan})
determine the Littlewood-Richardson coefficients as the number of
equivalence classes.
In
a similar way, fusion coefficients would count equivalence classes under
fusion-Knuth relations. This, however, remains to be the subject of
further
investigation.

\section*{Acknowledgments}  

I am indebted to Nantel Bergeron for constant guidance and
advice throughout this work and Terry Gannon who first suggested
the project and gave insights into fusion theory. I am also grateful
to Mark Shimozono and Stephanie van Willigenburg for helpful comments.

\end{document}